\documentclass{amsart}
\usepackage{amssymb,amsmath,amscd,enumerate,verbatim}
\usepackage[colorlinks=true,linkcolor=blue,citecolor=blue]{hyperref}
\usepackage{cleveref,graphicx}
\usepackage[all]{xy}

\numberwithin{equation}{section}

\newtheorem{thm}{\bf Theorem}[section]
\newtheorem{lem}[thm]{\bf Lemma}

\newtheorem{cor}[thm]{\bf Corollary}
\newtheorem{prop}[thm]{\bf Proposition}

\theoremstyle{definition}

\newtheorem{rem}[thm]{Remark}

\newtheorem*{thm*}{Theorem}
\newtheorem{lem*}[thm]{\bf Lemma}

\DeclareMathOperator{\Ass}{Ass}

\DeclareMathOperator{\Spec}{Spec}

\DeclareMathOperator{\Proj}{Proj}

\DeclareMathOperator{\Res}{Res}

\newcommand{\Z}{{\mathbb Z}}
\newcommand{\N}{{\mathbb N}}
\newcommand{\Pb}{\mathbb{P}}
\newcommand{\Q}{{\mathbb Q}}

\newcommand{\kk}{\Bbbk}

\def\pp{{\mathfrak p}}

\subjclass[2010]{13F20, 14N05, 13A02}

\begin{document}

\title[A sharp bound for the resurgence]{A sharp bound for the resurgence of \\  sums of ideals}
\author[D.V. Kien]{Do Van Kien}
\address{Department of Mathematics, Hanoi Pedagogical University 2, Vinh Phuc, Viet Nam}
\email{dovankien@hpu2.edu.vn}
\author[H.D. Nguyen]{Hop D. Nguyen}
\address{Institute of Mathematics, Vietnam Academy of Science and Technology, 18 Hoang Quoc Viet, 10307 Hanoi, Vietnam}
\email{ngdhop@gmail.com}
\author[L.M. Thuan]{Le Minh Thuan}
\address{Department of Mathematics, Hanoi Pedagogical University 2, Vinh Phuc, Viet Nam}
\email{leminhthuan1998cp@gmail.com}

\begin{abstract}
We prove a sharp upper bound for the resurgence of sums of ideals involving disjoint sets of variables, strengthening work of Bisui--H\`a--Jayanthan--Thomas. Complete solutions are delivered for two conjectures proposed by these authors. For given real numbers $a$ and $b$, we consider the set $\Res(a,b)$ of possible values of the resurgence of $I+J$ where $I$ and $J$ are ideals in disjoint sets of variables having resurgence $a$ and $b$, respectively. Some questions and partial results about $\Res(a,b)$ are discussed.
\end{abstract}

\maketitle

\section{Introduction}

Let $\kk$ be a field, $A=\kk[x_1,\ldots,x_d]$ be a polynomial ring with variables of degree 1. Let $I$ be an ideal of $A$. Denote by $I^{(n)}$ the $n$-th symbolic power of $I$ (defined in terms of associated primes):
\begin{equation}
\label{eq_ass_symbolpow}
I^{(n)}=\bigcap_{\pp \in \Ass(I)} \left(I^nA_\pp \cap A\right).
\end{equation}
A classical problem in commutative algebra is the comparison between ordinary and symbolic powers. One of the celebrated results in this area is
\begin{thm}[Hochster--Huneke, Ein--Lazarsfeld--Smith {\cite{ELS01, HH02}}]
\label{thm_HH}
Let $I$ be a homogeneous ideal in a polynomial ring $A=\kk[x_1,\ldots,x_d]$, where $d\ge 2$. Denote by $h$ the big height of $I$, namely the maximal height of an associated prime of $I$. Then for all $n\ge 1$, there is a containment $I^{(hn)}\subseteq I^n$. In particular, we always have $I^{((d-1)n)}\subseteq I^n$ for all $n\ge 1$.
\end{thm}

The resurgence and asymptotic resurgence of $I$, denoted by $\rho(I)$ and $\rho_a(I)$, were introduced in \cite{BH10, GHV13} as follows:
\begin{align*}
	\rho(I)  &= \sup\left\{\frac{m}{r}: I^{(m)} \not\subseteq I^r, m\ge 1,r\ge 1 \right \},\\
	\rho_a(I) &= \sup\left\{\frac{m}{r}: I^{(mt)} \not\subseteq I^{rt}  ~\text{for $t\gg 0$} \right\}.
\end{align*}
 These numbers measure the (asymptotic) difference between the ordinary and symbolic powers of $I$. Recently, many authors have paid attention to the resurgence and asymptotic resurgence  of homogeneous ideals; see, for example \cite{BHJT, BH10, CHHVT, DFMS, DHNSTT, GHM2020,GHV13, JKM2021}. We refer to \cite{DDS+} for a survey on various topics concerning symbolic powers.

It is clear that $1\le \rho_a(I) \le \rho(I)$. Some additional bounds for $\rho(I)$ and $\rho_a(I)$ are given in \cite{BH10,DFMS,GHV13, HTr19, JKM2021}, for instance if $I$ is a squarefree monomial ideal, then $\rho(I)\le d(I)$, the maximal generating degree of $I$ by \cite[Theorem 3.18]{DFMS} and \cite[Corollary 3.6]{HTr19}. As far as we know, there is still no effective algorithm to compute $\rho(I)$ and $\rho_a(I)$ for a general ideal $I$. Some algorithms for computing $\rho_a(I)$ when $I$ is squarefree monomial are available thanks to \cite[Section 2]{DFMS}, \cite[Theorem 5.3, Appendix A]{GSV2022} and \cite[Theorem 3.7]{Vil2022}. Apparently no analogous algorithm is known for the resurgence of an arbitrary monomial ideal.

Let $A, B$ be standard graded polynomial rings over $\kk$. Let $I\subseteq A$, $J\subseteq B$ be non-zero proper homogeneous ideals. In \cite{BHJT}, the resurgence and asymptotic resurgence of the ideal of $R=A\otimes_\kk B$ generated by $I$ and $J$, simply written as $I+J$, was studied. The construction of the ideal $I+J$ is a classical construction and it is relevant in commutative algebra and algebraic geometry since it corresponds to the notions of tensor products of $\kk$-algebras and fiber products of schemes over $\Spec(\kk)$. For various classical invariants, the value at $I+J$ (or $R/(I+J)$) are determined by the values at $I$ and $J$ (or the corresponding quotient rings). This is the case, for example, for the Krull dimension, depth, graded Betti numbers, and hence the Castelnuovo--Mumford regularity.  It turns out that we can also determine the asymptotic resurgence of $I+J$ from those of $I$ and $J$: It is proved in  \cite[Theorem 2.6]{BHJT} that there is an equality 
\begin{equation}
\label{eq_rhoa_sum}
\rho_a(I+J)=\max\{\rho_a(I),\rho_a(J)\}.
\end{equation}
On the other hand, for the resurgence of $I+J$, the best known information is given by the following inequalities \cite[Theorem 2.7]{BHJT}:
\begin{equation}
\label{eq_upperbound}
\max \{\rho(I),\rho(J)\} \le \rho(I+J) \le \rho(I)+\rho(J). 
\end{equation}
It is also noted in \cite[Remark 2.9]{BHJT} that the authors are not aware of any case where the upper bound is attained. Some partial improvements of \eqref{eq_upperbound} was provided in recent work by Jayanthan, Kumar, and Mukundan \cite[Theorems 3.6, 3.9]{JKM2021}. As our first main result, we prove the following improved upper bound for the resurgence of sums of ideals. Our result also confirms that the upper bound \eqref{eq_upperbound} is really strict: As $\rho(I), \rho(J)\ge 1$, the result below implies that $\rho(I+J) \le \rho(I)+\rho(J)-2/3$ always holds.
\begin{thm}[= \Cref{thm_sharpbound_rho_sum}]
Let $A, B$ be standard graded polynomial rings over $\kk$. Let $I\subseteq A$, $J\subseteq B$ be non-zero proper homogeneous ideals. Then there are inequalities
$$
\max \{\rho(I),\rho(J)\} \le \rho(I+J) \le \max \left\{\rho(I),\rho(J),\dfrac{2(\rho(I)+\rho(J))}{3}\right\},
$$
and the upper bound is sharp. In particular, if the inequality $\max\{\rho(I),\rho(J)\}\ge 2\min\{\rho(I),\rho(J)\}$ holds then there is an equality 
$\rho(I+J)=\max\{\rho(I),\rho(J)\}.$
\end{thm}
The last assertion of this result seems to be rather unexpected. The proof of \Cref{thm_sharpbound_rho_sum} is somewhat similar to, but differs in a crucial way, from the proof method of \cite[Theorem 2.7]{BHJT}, namely we employ more efficiently the binomial expansion formula for ``associated'' symbolic powers of $I+J$. This formula was proved first for ``minimal'' symbolic powers (defined in terms of minimal primes) in \cite[Theorem 3.4]{HNTT}, and later for ``associated'' symbolic powers as well in \cite[Theorem 4.1]{HJKN}. We stress that, in this paper, we focus solely on the ``associated'' symbolic powers defined by Formula \eqref{eq_ass_symbolpow}. For related results on the behaviour of (asymptotic) resurgence under taking sum, product, and intersection, we refer to  \cite{JKM2021}.

The next two main results answer completely two conjectures proposed by Bisui--H\`a--Jayanthan--Thomas in \cite[Conjectures 3.8 and 3.9]{BHJT}, one negatively and the other positively. \Cref{thm_main} yields the following negative answer to Conjecture 3.8 in \emph{ibid.}
\begin{cor}[= \Cref{cor_boundedrho}]
 Let $I$ be a nonzero proper homogeneous ideal in $A=\kk[x_1,\ldots,x_d]$. Let $X=\Proj(A/I)$ and for each $m\ge 1$, let $I^{[m]}$ be the defining ideal of the fiber product $\underbrace{X\times_\kk \cdots \times_\kk X}_{\textup{$m$ times}}$, as a closed subscheme of $\underbrace{\Pb^{d-1}_\kk \times_\kk \cdots \times_\kk \Pb^{d-1}_\kk}_{\textup{$m$ times}}$. Then for all $m\ge 1$, the strict inequality $\rho(I^{[m]})<2\rho(I)$ always holds.
\end{cor}

There are equalities $0\le \rho(I)-\rho_a(I)\le \dim(A)-1$ thanks to \Cref{thm_HH}. But so far it is not clear whether by varying the number of variables, $\rho(I)-\rho_a(I)$ may become arbitrarily large.  This issue was raised in \cite[Conjecture 3.9]{BHJT}. Using a formula on the asymptotic resurgence of the monomial ideal of certain star configurations in \cite[Theorem C]{LBG15}, we give a positive answer to this conjecture. We show that for certain sequence $(P_m)$ of squarefree monomial ideals, each generated in a single degree, the difference $\rho(P_m)-\rho_a(P_m)$ tends to infinity.
\begin{thm}[Cf. \Cref{thm_large_difference}]
\label{thm_main}
Let $m\ge 2$ be an integer, $R=\kk[x_{1,i}, x_{2,i}, x_{3,i}: 1\le i\le 2m-1]$ and for each $1\le j\le 3$, let $I_j$ be the ideal generated by all the products of $m$ distinct variables among the $(2m-1)$ variables $x_{j,1},\ldots,x_{j,2m-1}$:
\[
I_j=(x_{j,i_1}x_{j,i_2}\cdots x_{j,i_m}: 1\le i_1 <i_2<\cdots<i_m \le 2m-1).
\]
Denote $P_m=I_1+I_2+I_3$. Then $\rho(P_m)\ge \dfrac{3m}{4}$ and  $\rho_a(P_m)=\dfrac{m^2}{2m-1}$. In particular, $\dfrac{2m^2-3m}{4(2m-1)} \le \rho(P_m)-\rho_a(P_m)$, so $\lim\limits_{m\to \infty}\left(\rho(P_m)-\rho_a(P_m)\right)= \infty$.
\end{thm}
\textbf{Organization.} The new upper bound for the resurgence of sums of ideals is given in \Cref{sect_sums}. We construct for any given integer $d\ge 1$, a monomial ideal $I$ in three variables with resurgence number $1$, such that the equality between ordinary and symbolic powers holds for the first $d$ powers, but fails for the $(d+1)$ ones (\Cref{lem_rho1_3generators}). This seemingly new construction is useful for showing that the upper bound of \Cref{thm_main} is sharp even in the non-trivial case $\max\{\rho(I),\rho(J)\}< 2\min\{\rho(I),\rho(J)\}$, e.g. when $\rho(I)=\rho(J)$. As an application of the sharp bound for the resurgence of sums, we deduce \Cref{cor_boundedrho}. That the difference $\rho(I)-\rho_a(I)$ can be arbitrarily large is proved in \Cref{sect_difference}. In the last \Cref{sect_resset}, for given real numbers $a$ and $b$, we study the set 
\begin{gather*}
\Res(a,b)=\left\{\rho(I+J):  \text{$I$ and $J$ live  in disjoint sets of variables} \right.\\
                                        \qquad   \text{such that} \left. \rho(I)=a, \rho(J)=b \right\}. 
\end{gather*}
Building upon a previous work \cite{JKM2021} and \Cref{lem_rho1_3generators}, we describe completely $\Res(1,1)$. We end the paper by discussing some unsolved questions on the set $\Res(a,b)$.

\section{Sharp upper bound for the resurgence of sums of ideals}
\label{sect_sums}
Let $A$ be a noetherian ring and $I$ an ideal of $A$. In this paper, we only work with the following notion of the $n$-th symbolic power of $I$:
\[
I^{(n)}=\bigcap_{\pp \in \Ass(I)} (I^nA_\pp \cap A).
\]
Thus symbolic powers in our sense are defined in terms of \emph{associated primes}, not minimal primes. It is well-known that $I^{(n)} \subseteq I^{(1)}=I$ always holds for every $n\ge 1$.

The following result is from recent work by H\`a et al.
\begin{thm}[{\cite[Theorem 4.1]{HJKN}}]
	\label{thm.symbBIN} Let $A, B$ be standard graded polynomial rings over a field $\kk$ and $I \subseteq A, J \subseteq B$ be nonzero proper homogeneous ideals. Then, for any $s \in \N$, we have
	\begin{align} \label{eq.symbBIN}
		(I+J)^{(s)} = \sum_{i=0}^s I^{(i)}J^{(s-i)}.
	\end{align}
\end{thm}

\begin{lem}
\label{lem_sum_and_resurgence}
Let $A_1,\ldots,A_p$ be standard graded polynomial rings over $\kk$, and $I_i\subseteq A_i$ be a proper homogeneous ideal. Assume that for $1\le i\le p$, $m_i,r_i\ge 1$ be integers such that $I_i^{(m_i)}\not\subseteq I_i^{r_i}$. Denote $P=I_1+\cdots+I_p\subseteq A_1\otimes_\kk \cdots \otimes_\kk A_p$. Then
\[
P^{(m_1+\cdots+m_p)} \not\subseteq P^{r_1+\cdots+r_p-p+1}.
\]
In particular, if $p=v+1$, $m_1=\cdots=m_{v+1}=m$, and $r_1=\cdots=r_{v+1}=v$ for some $v\ge 1$, then $P^{(m(v+1))} \not\subseteq P^{v^2}$.
\end{lem}
\begin{proof}
For the first assertion, it suffices to consider the case $p=2$, as the general case follows by induction on $p$. In this case, denote $A=A_1,I=I_1,B=A_2,J=I_2$ for simplicity, so $P=I+J\subseteq R=A\otimes_\kk B$. 

We have to show that $P^{(m_1+m_2)}\not\subseteq P^{r_1+r_2-1}$ given that $I^{(m_1)} \not\subseteq I^{r_1}$ and $J^{(m_2)}\not\subseteq J^{r_2}$. From  \cite[Lemma 3.2]{BHJT}, as soon as $f\in A, g\in B$, $f\in I^{(m_1)} \setminus I^{r_1}$, and $g\in J^{(m_2)}\setminus J^{r_2}$, 
\[
fg \notin P^{r_1+r_2-1}.
\]
Clearly $fg \in I^{(m_1)}J^{(m_2)} \subseteq P^{(m_1+m_2)}$, so $P^{(m_1+m_2)}\not\subseteq P^{r_1+r_2-1}$. The first assertion follows. The remaining assertion is a simple accounting.
\end{proof}
The result \cite[Theorem 2.7]{BHJT} has shown that $\rho(I+J)\le \rho(I)+\rho(J)$. The following statement gives a better and sharp upper bound.
 
\begin{thm}
\label{thm_sharpbound_rho_sum}
Let $A, B$ be standard graded polynomial rings over $\kk$. Let $I\subseteq A$, $J\subseteq B$ be non-zero proper homogeneous ideals. Then there are inequalities
\[
\max \{\rho(I),\rho(J)\} \le \rho(I+J) \le \mathop{\sup_{m, n\in \Z}}_{m,n\ge 2}\left\{\rho(I),\rho(J), \dfrac{m\rho(I)+n\rho(J)}{m+n-1} \right\}.
\]
More precisely, this is equivalent to 
$$
\max \{\rho(I),\rho(J)\} \le \rho(I+J) \le \max \left\{\rho(I),\rho(J),\dfrac{2(\rho(I)+\rho(J))}{3}\right\},
$$
and the upper bound is sharp. Furthermore, if $\max\{\rho(I),\rho(J)\}\ge 2\min\{\rho(I),\rho(J)\}$ then $\rho(I+J)=\max\{\rho(I),\rho(J)\}$.
\end{thm}
The inequality $\max \{\rho(I),\rho(J)\} \le \rho(I+J)$ was proved in \cite[Theorem 2.7]{BHJT}.  We remark that the upper bound in Theorem \ref{thm_sharpbound_rho_sum} is better than that one in \cite[Theorem 2.7]{BHJT}. In fact, as $\min\{\rho(I),\rho(J)\}\ge 1$, 
\begin{align*}
\max\left\{\rho(I), \rho(J), \frac{2}{3}(\rho(I)+\rho(J)) \right\}&\le \rho(I)+\rho(J)-\dfrac{2}{3}\min\left\{\rho(I), \rho(J)\right\} \\
                                                                  &\le \rho(I)+\rho(J)-2/3. 
\end{align*}
\Cref{thm_sharpbound_rho_sum} therefore gives an explanation to \cite[Remark 2.9]{BHJT}. 

The proof of Theorem \ref{thm_sharpbound_rho_sum} employs the following lemma.
\begin{lem}
\label{lem_evaluate_max}
Let $a, b$ be non-negative real numbers. Then 
\[
\mathop{\sup_{m, n\in \Z}}_{m,n\ge 2}\left\{a,b, \dfrac{ma+nb}{m+n-1} \right\}=\max \left\{a,b,\dfrac{2(a+b)}{3} \right\}.
\]
\end{lem}
\begin{proof}
Choosing $m=n=2$, we see that the left-hand side is not smaller than the right-hand side. Hence it remains to prove the reverse inequality.

It is harmless to assume that $a\ge b$. Note that 
\[
\max \left\{a,b,\dfrac{2(a+b)}{3} \right\} =\begin{cases}
                                             a, &\text{if $a\ge 2b$},\\
                                             2(a+b)/3, &\text{if $b\le a < 2b$}.
          
          \end{cases}
\]
For integers $m,n\ge 2$, since $a-b,b\ge 0$, 
\begin{align*}
\dfrac{ma+nb}{m+n-1} & = \dfrac{ma-(m-1)b}{m+n-1}+b  =\dfrac{m(a-b)+b}{m+n-1}+b  \\
                     & \le \dfrac{m(a-b)+b}{m+1}+b = a + \dfrac{2b-a}{m+1}.
\end{align*}
If $a\ge 2b$, then $a + \dfrac{2b-a}{m+1} \le a.$ If $a<2b$, using $m\ge 2$,
\[
 a + \dfrac{2b-a}{m+1} \le a + \dfrac{2b-a}{3}=\dfrac{2(a+b)}{3}.
\]
This finishes the proof.
\end{proof}

\begin{proof}[Proof of Theorem \ref{thm_sharpbound_rho_sum}]
Take $h,r\ge 1$ such that
\[
h/r>D=\mathop{\sup_{m, n\in \Z}}_{m,n\ge 2}\left\{\rho(I),\rho(J), \dfrac{m\rho(I)+n\rho(J)}{m+n-1} \right\}.
\]
We show that $(I+J)^{(h)} \subseteq (I+J)^r$. Note that $(I+J)^{(h)} \subseteq (I+J)^{(1)}=I+J$, so it suffices to consider the case $r\ge 2$. Without loss of generality, we assume $\rho(I)\ge \rho(J)$.

By \Cref{thm.symbBIN}, there is an equality $(I+J)^{(h)}=\sum_{i=0}^h I^{(i)}J^{(h-i)}$. Using this, we want to show that $I^{(i)}J^{(h-i)}\subseteq (I+J)^r$ for all $0\le i\le h$. If $i=0$, $h/r> \rho(I)$, so $I^{(h)}\subseteq I^r$. Similarly, if $i=h$ then $J^{(h)}\subseteq J^r$. Hence it remains to consider the case $1\le i\le h-1$.

There exists a unique integer $m\ge 0$ such that $m\rho(I) < i \le (m+1)\rho(I)$. Note that $I^{(i)}\subseteq I^m$. If $m\ge r-1$, then using $h-i\ge 1$,
$$
I^{(i)}J^{(h-i)}\subseteq I^mJ\subseteq I^{r-1}J\subseteq (I+J)^r,
$$
so we are done. Assume that $m\le r-2$.

If $m=0$, then $i\le \rho(I)$. Since $h/r>\rho(I)$, we get $h-i\ge h-\rho(I)>(r-1)\rho(I)\ge (r-1)\rho(J)$. Therefore $I^{(i)}J^{(h-i)}\subseteq I J^{r-1} \subseteq (I+J)^r$.

Assume that $1\le m\le r-2$. Denote $n=r-m$, then $n\ge 2$ and $r=(m+1)+n-1$. Since 
$$
h/r> D\ge  \dfrac{(m+1)\rho(I)+n\rho(J)}{m+n}=  \dfrac{(m+1)\rho(I)+n\rho(J)}{r},
$$
it follows that $h>(m+1)\rho(I)+n\rho(J)$. Thus $h-i\ge h-(m+1)\rho(I)> n\rho(J)$. This yields $J^{(h-i)}\subseteq J^n$, consequently $I^{(i)}J^{(h-i)} \subseteq I^mJ^n \subseteq (I+J)^r$. Hence $\rho(I+J)\le D$.

The second assertion holds since by Lemma \ref{lem_evaluate_max}, 
$$
D=\max \left\{\rho(I),\rho(J),\dfrac{2(\rho(I)+\rho(J))}{3}\right\}.
$$
That the upper bound is sharp follows from part (2) of Lemma \ref{lem_example_sharp} below, where we give an example with $\rho(I)=\rho(J)=1$ and $\rho(I+J)=4/3$.

When $\max\{\rho(I),\rho(J)\}\ge 2\min\{\rho(I),\rho(J)\}$, we have
\[
\max \left\{\rho(I),\rho(J),\dfrac{2(\rho(I)+\rho(J))}{3}\right\}=\max\{\rho(I),\rho(J)\}. 
\]
Hence $\rho(I+J)=\max \left\{\rho(I),\rho(J)\right\}$ in this case. The proof is completed.
\end{proof}

\begin{lem}
\label{lem_example_sharp}
Let $A=\kk[x,y,z]$, $I=(x^3,xy^2,y^3)\cap (x,z)=(x^3,xy^2,y^3z)$. 
\begin{enumerate}[\quad \rm (1)]
 \item For all $n\ge 2$, there is a chain $I^{(n)}=I^n+x^3y^3I^{n-2}\subseteq I^{n-1}$. In particular, $\rho(I)=1$. Moreover $x^3y^3\in I^{(2)}\setminus I^2$.

\item Let  $B=\kk[u,v,w]$, and  $J=(u^3,uv^2,v^3w)\subseteq B$. Then $\rho(J)=\rho(I)=1$ while $\rho(I+J)= 4/3$.
\end{enumerate}

\end{lem}
\begin{proof}
(1): This follows from the more general \Cref{lem_rho1_3generators} below.

(2): By part (1), $x^3y^3 \in I^{(2)}\setminus I^2, u^3v^3\in J^{(2)}\setminus J^2$.  This implies 
\[
x^3y^3u^3v^3\in I^{(2)}J^{(2)}\setminus (I^2+J^2) \subseteq (I+J)^{(4)}\setminus (I+J)^3.
\]
Hence $\rho(I+J)\ge 4/3$ and by \Cref{thm_sharpbound_rho_sum}, $\rho(I+J)=4/3$.
\end{proof}
\begin{lem}
\label{lem_rho1_3generators}
Let $d\ge 1$ be an integer, $I=(x^{2d+1},x^{2d-1}y^2,y^{2d+1}z)\subseteq A=\kk[x,y,z]$. Then the following statements hold.
\begin{enumerate}[\quad \rm (1)]
 \item For each $1\le n\le d$, there is an equality $I^{(n)}=I^n$.

 \item For each $n\ge d+1$, there is a chain 
 $$
 I^{(n)}=I^n+(x^dy)^{2d+1} I^{n-d-1} \subseteq I^{n-1}.
 $$
 In particular, $\rho(I)=1$. Moreover $(x^dy)^{2d+1} \in I^{(d+1)}\setminus I^{d+1}$.
 \end{enumerate}
\end{lem}

\begin{proof}
The irredundant primary decomposition 
$$I=(x^{2d+1},x^{2d-1}y^2,y^{2d+1}) \cap (x^{2d-1},z)
$$ 
implies that for all $n\ge 1$,
\begin{equation}
 \label{eq_symbpow_formula}
 I^{(n)}=(x^{2d+1},x^{2d-1}y^2,y^{2d+1})^n \cap (x^{2d-1},z)^n.
\end{equation}

\textsc{Step 1:} Let $f=x^ay^bz^c$ be a monomial in $I^{(n)}$. Note that belonging to the ideal $(x^{2d+1},x^{2d-1}y^2,y^{2d+1})^n$, $x^ay^b$ has a divisor $(x^{2d+1})^g(x^{2d-1}y^2)^h(y^{2d+1})^i=x^{(2d+1)g+(2d-1)h}y^{2h+(2d+1)i}$, where 
\begin{equation}
\label{eq_hypothesis0}
g,h,i\ge 0, g+h+i= n.
\end{equation}
We deduce
\begin{align}
a & \ge (2d+1)g+(2d-1)h, \label{eq_hypothesis1}\\
b &\ge 2h+(2d+1)i, \label{eq_hypothesis2}  \\ 
\left\lfloor \frac{a}{2d-1} \right \rfloor+c & \ge n. \label{eq_hypothesis3}
\end{align}
The last inequality holds since $f\in (x^{2d-1},z)^n$.
 
Adding \eqref{eq_hypothesis1} and \eqref{eq_hypothesis2}, then using \eqref{eq_hypothesis0}, it follows that
\begin{equation}
\label{eq_alarge}
a \ge (2d+1)n-b. 
\end{equation}
Note that $x^ay^bz^c \in I^n$ if and only if it has a divisor $(x^{2d+1})^p(x^{2d-1}y^2)^q(y^{2d+1}z)^r$ where $p,q,r\ge 0, p+q+r=n$. Equivalently,  $f\in I^n$ if and only if the following system has an integral solution
\[
\begin{cases}
p,q,r &\ge 0, p+q+r=n,\\
a &\ge (2d+1)p+(2d-1)q,\\
 b &\ge 2q+(2d+1)r,\\
 c& \ge r.
\end{cases}
\]
Since $(2d+1)p+(2d-1)q=(2d+1)n-(2q+(2d+1)r)$ and $r=n-(p+q)$, the last system has an integral solution $(p,q,r)$ if and only if the following system has an integral solution $(p,q)$: 
\begin{align}
p,q &\ge 0, \label{eq_system1}\\
n-c &\le p+q \le n, \label{eq_system2}\\
(2d+1)n-b &\le (2d+1)p+(2d-1)q \le a. \label{eq_system3}
\end{align}

\begin{figure}[ht!]
\includegraphics[width=78ex]{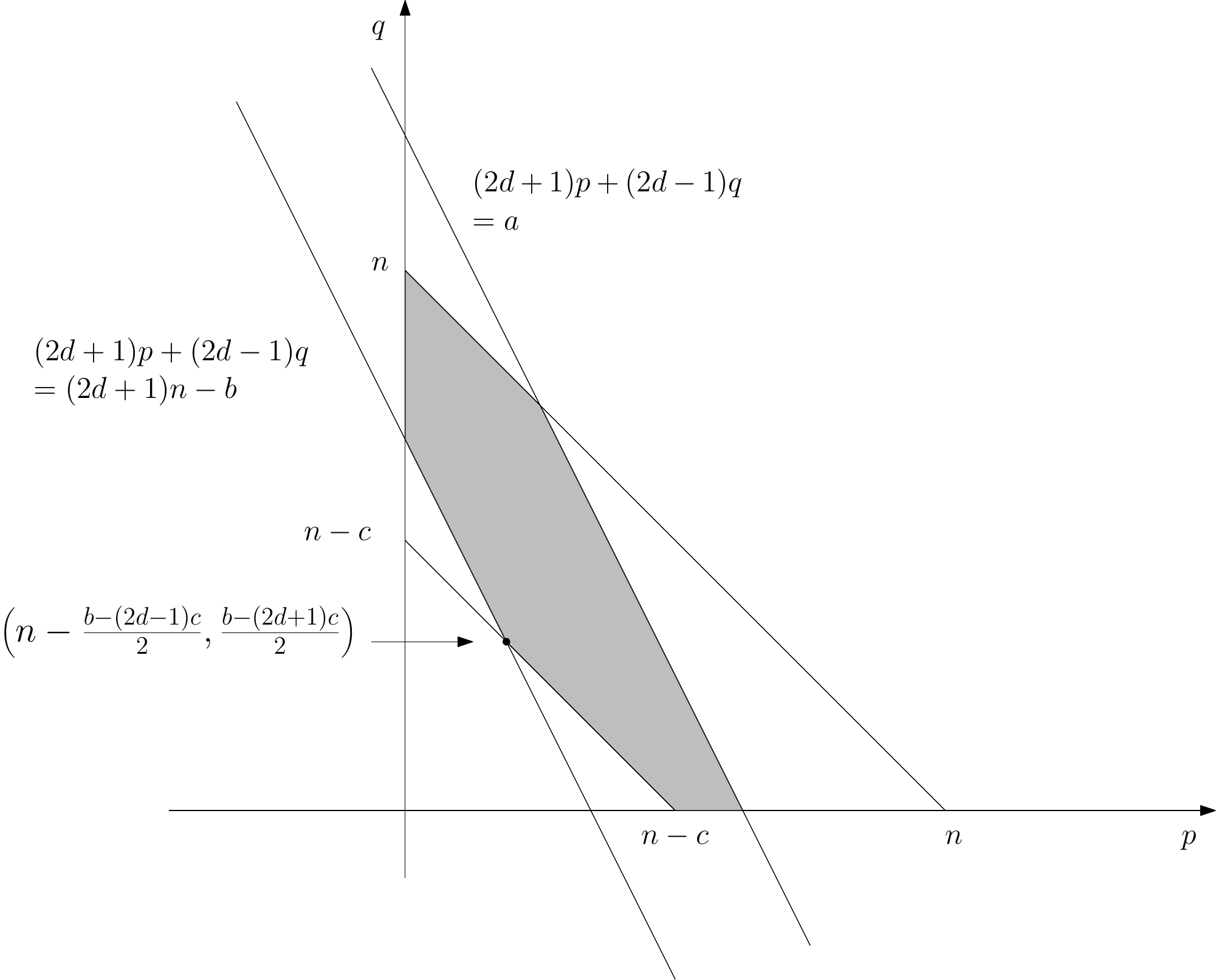}
\caption{The representation of the system \eqref{eq_system1}--\eqref{eq_system3} in the coordinate plane.}
\label{fig_inequal}
\end{figure}

\textsc{Step 2}: We have a crucial observation.

\textbf{Claim 1:} For any $n\ge 1$, if $f\notin I^n$ then the following conditions are simultaneously satisfied
\begin{equation}
 \label{eq_condition_nosol}
 i\ge c+1, \, \text{$b\le 2n+(2d-1)c-1$ and $b+c$ is odd}.
\end{equation}
\emph{Proof of Claim 1}: If $i\le c$, then $g+h=n-i\ge n-c$. Thus the system \eqref{eq_system1}--\eqref{eq_system3} has a solution $(p,q)=(g,h)$, thanks to the hypotheses \eqref{eq_hypothesis0},   \eqref{eq_hypothesis1}, \eqref{eq_hypothesis2}. It remains to consider the case $i\ge c+1$. 

Note that this yields $c\le i-1\le n-1$ and \eqref{eq_hypothesis2} implies that
\begin{equation}
\label{eq_blarge}
b\ge 2h+(2d+1)i\ge 2h+(2d+1)(c+1) \ge 2d+1.
\end{equation}

\textsf{Case 1:} $b\ge 2n+(2d-1)c$. The system \eqref{eq_system1}--\eqref{eq_system3} has a solution $(p,q)=(0,n-c)$, as per \eqref{eq_hypothesis3}, $(2d-1)(n-c)\le a$.

\textsf{Case 2:} $b\le 2n+(2d-1)c-1$ and $b+c$ is even. Solving for
\[
\begin{cases}
 p+q&= n-c,\\
 (2d+1)p+(2d-1)q &=(2d+1)n-b,
\end{cases}
\]
we get
\[
(p,q)=\left(n-\dfrac{b-(2d-1)c}{2},\dfrac{b-(2d+1)c}{2}\right).
\]
Since $b+c$ is even, $p,q\in \Z$. Since $b\le 2n+(2d-1)c-1$, $p$ is non-negative, and thanks to \eqref{eq_blarge}, it follows that $q\ge 0$. Hence \eqref{eq_system1}--\eqref{eq_system3} admits an integral solution.

Now we are left with the case $i\ge c+1$, $b\le 2n+(2d-1)c-1$ and $b+c$ is odd, namely Claim 1 is true.

\textsc{Step 3:} Assume that $I^{(n)}\neq I^n$ for some $n \ge 1$. There exists a monomial $f=x^ay^bz^c \in I^{(n)}\setminus I^n$. Choose the integers $g,h,i$ as in Step 1. By Claim 1, we have $i\ge c+1$ and $b\le 2n+(2d-1)c-1$.

Using \eqref{eq_hypothesis2},
\[
2n+(2d-1)c-1 \ge b \ge 2h+(2d+1)i \ge 2h+(2d+1)(c+1).
\]
Simplifying, this yields,
\begin{equation}
\label{eq_nlarge}
n\ge h+c+d+1 \ge d+1.
\end{equation}
In particular, this shows that $I^{(n)}=I^n$ for all $1\le n\le d$.

\textsc{Step 4}: We prove for each $n\ge d+1$ that
\[
I^{(n)}=I^n+(x^dy)^{2d+1}I^{n-d-1}.
\]
The containment $I^n+(x^dy)^{2d+1}I^{n-d-1} \subseteq I^{(n)}$ is elementary. Indeed, it suffices to show that $f_0=(x^dy)^{2d+1}=x^{2d^2+d}y^{2d+1} \in I^{(d+1)}$. In fact $f_0=(x^{2d+1})^dy^{2d+1} \in (x^{2d+1},x^{2d-1}y^2,y^{2d+1})^{d+1}$ and 
$$
f_0\in (x^{2d^2+d-1})  \subseteq (x^{2d-1},z)^{d+1}.
$$
Hence $f_0 \in I^{(d+1)}$.

For the reverse containment, take any monomial $f=x^ay^bz^c \in I^{(n)}$. Define the numbers $g,h,i$ as in Step 1.  Assuming $f\notin I^n$, we claim that $f\in (x^dy)^{2d+1}I^{n-d-1}$. By Claim 1, $i\ge c+1$, $b\le 2n+(2d-1)c-1$ and $b+c$ is odd.

\textbf{Claim 2:} We have inequalities $b\ge 2d+1, a\ge 2d^2+d$.

That $b\ge 2d+1$ follows from \eqref{eq_blarge}. Assume that $a\le 2d^2+d-1=(2d-1)(d+1)$, then $a/(2d-1)\le d+1$. As in \eqref{eq_nlarge}, we get $n-c\ge h+d+1$. Hence together with \eqref{eq_hypothesis2},
\[
h+d+1 \le n-c \le \dfrac{a}{2d-1} \le d+1.
\]
This implies that $h=0$ and $d+1=n-c=\dfrac{a}{2d-1}$. Again by \eqref{eq_hypothesis2},
\[
(2d+1)i\le b \le 2n+(2d-1)c-1 = 2n+(2d-1)(n-d-1)-1=(2d+1)(n-d),
\]
so $i\le n-d$. But then $g=n-h-i\ge d$, and \eqref{eq_hypothesis1} yields
\[
a\ge (2d+1)g \ge 2d^2+d,
\]
a contradiction. Thus the above assumption is wrong, and $a\ge 2d^2+d$.

Now $f=x^ay^bz^c=x^{2d^2+d}y^{2d+1} x^{a'}y^{b'}z^c$, where $a'=a-(2d^2+d), b'=b-(2d+1)$. Denote $n'=n-d-1$, then \eqref{eq_nlarge} implies that $n'\ge c\ge 0$. We wish to show that $x^{a'}y^{b'}z^c \in I^{n'}$. As in Step 1, this means the following system has an integral solution
\[
\begin{cases}
p',q' &\ge 0,\\
n'-c &\le p'+q'\le n',\\
(2d+1)n'-b' &\le (2d+1)p'+(2d-1)q'\le a'.
\end{cases}
\]
Solving for
\[
\begin{cases}
p'+q' &= n'-c,\\
(2d+1)p'+(2d-1)q'&=(2d+1)n'-b',
\end{cases}
\]
we get that
\[
q'=\dfrac{b'-(2d+1)c}{2}=\dfrac{b-(2d+1)(c+1)}{2},\quad p'=n'-c-q'. 
\]
Since $b+c$ is odd, $p', q'\in \Z$. It remains to check that $0\le q'\le n'-c$. The first inequality holds thanks to \eqref{eq_blarge}. The second inequality can be rewritten as
\[
b-(2d+1)(c+1)\le 2(n-d-1-c),
\]
equivalently $b \le 2n+(2d-1)c-1$, which is valid by \eqref{eq_condition_nosol}. Hence the desired containment $f\in (x^dy)^{2d+1}I^{n-d-1}$ holds, and we finish the proof that for every $n\ge d+1$,
\[
I^{(n)}=I^n+(x^dy)^{2d+1}I^{n-d-1}.
\]
\textsc{Step 5:} To prove (2), it remains to show that for each $n\ge d+1$,
\[
I^n+(x^dy)^{2d+1}I^{n-d-1} \subseteq I^{n-1}.
\]
This is clear since $(x^dy)^{2d+1} \in ((x^{2d+1})^d) \subseteq I^d$. Now $\rho(I)=1$ follows since $I^{(n)}\subseteq I^{n-1}$ for all $n\ge 1$. Finally, we check that $(x^dy)^{2d+1} \in I^{(d+1)}\setminus I^{d+1}$. The containment was established in Step 4. Assume that $(x^dy)^{2d+1}\in I^{d+1}=(x^{2d+1},x^{2d-1}y^2,y^{2d+1}z)^{d+1}$. Inspecting supports, we deduce
\[
(x^dy)^{2d+1} \in (x^{2d+1},x^{2d-1}y^2)^{d+1},
\]
and after simplifying common factors,
\[
xy^{2d+1} \in  (x^2,y^2)^{d+1}.
\]
This is clearly a contradiction. The proof is completed.
\end{proof}

\begin{rem}
\label{rem_example_strict}
The two inequalities in \Cref{thm_sharpbound_rho_sum} can be both strict. For example, let $I=(x^3,xy^2,y^3z)\subseteq A=\kk[x,y,z]$ and 
$$
J=(t^5,t^3u^2,u^5v)\subseteq B=\kk[t,u,v].
$$
By \Cref{lem_rho1_3generators} for $d=2$,  $\rho(J)=1$ and $t^{10}u^5\in J^{(3)}\setminus J^3$.

We have two strict inequalities in the chain
$$
\max\{\rho(I),\rho(J)\}=1 < \rho(I+J)=\dfrac{5}{4}< \max \left\{\rho(I),\rho(J),\dfrac{2(\rho(I)+\rho(J))}{3}\right\}=\dfrac{4}{3}.
$$
Indeed, this can be seen using \cite[Theorem 3.9]{JKM2021}, or by direct arguments as follows. Denote $P=I+J$. We have $x^3y^3\in I^{(2)}\setminus I^2$, $t^{10}u^5\in J^{(3)}\setminus J^3$, hence 
\[
x^3y^3t^{10}u^5 \in I^{(2)}J^{(3)} \setminus (I^2+J^3) \subseteq P^{(5)}\setminus P^4.
\]
Hence $\rho(P)\ge 5/4$. Using the fact that $I^{(n+1)}\subseteq I^n, J^{(n+1)}\subseteq J^n$ for all $n\ge 1$, we get that $P^{(n)}\subseteq P^{n-2}$ for all $n\ge 2$. Moreover by direct inspection $P^{(n)}\subseteq P^{n-1}$ for $1\le n\le 4$. Hence 
\[
\rho(P)\le \mathop{\sup_{2\le m}}_{n\ge 5} \left\{\frac{m}{m}, \frac{n}{n-1}\right\} =\frac{5}{4}.
\]
The desired conclusion follows.
\end{rem}

Now we give a complete answer to \cite[Conjecture 3.8]{BHJT} about the resurgence numbers of iterated sums of an ideal. For an ideal $I$ in a polynomial ring $A$, the \emph{$d$-th iterated sum} $I^{[d]}$ of $I$ is the defining ideal of the tensor product 
$$
\underbrace{(A/I)\otimes_\kk (A/I) \otimes_\kk \cdots \otimes_\kk (A/I)}_{\text{$d$ times}}.
$$
as a quotient ring of $A^{\otimes d}$. For example, for $I = (x^2, xy) \subseteq \kk[x,y]$, the first three iterated sums of $I$ are given by:
\begin{itemize}
	\item  $I^{\left[1\right]}= ( x_1^2, x_1y_1 ) \subseteq \kk[x_1, y_1]$,
	\item  $I^{\left[2\right]}= ( x_1^2, x_1y_1, x_2^2, x_2 y_2 ) \subseteq \kk[x_1, y_1, x_2, y_2]$,
	\item  $I^{\left[3\right]}= ( x_1^2, x_1y_1,  x_2^2, x_2y_2,  x_3^2, x_3y_3 ) \subseteq \kk[x_1, y_1, x_2, y_2, x_3, y_3]$.
\end{itemize}
From \Cref{thm_sharpbound_rho_sum}, \cite[Conjecture 3.8]{BHJT} admits a negative answer.
\begin{cor}
\label{cor_boundedrho}
Let $I$ be any nonzero proper homogeneous ideal in a polynomial ring $A$. Then for all $d \in \mathbb{N}$ we have a strict inequality
$$\rho(I^{\left[d\right]}) < 2\rho(I).$$
\end{cor}
\begin{proof}
 We will prove the above inequality by induction. By definition of iterated sums, we have
$$\rho(I^{\left[d+1\right]})=\rho(I + I^{\left[d\right]}) \leq \max \left\lbrace \rho(I), \rho(I^{\left[d\right]}),\dfrac{2}{3}(\rho(I)+ \rho(I^{\left[d\right]})) \right \rbrace$$
for all $d \in \mathbb{N}, d \geq 1$. For $d=1$,
$$\rho(I^{\left[2\right]}) \leq \max \{\rho(I^{\left[1\right]}), \rho(I^{\left[1\right]}), \dfrac{4}{3}\rho(I^{\left[1\right]})\}=\dfrac{4}{3}\rho(I^{\left[1\right]}) < 2\rho(I).$$
Assume that the conclusion is true up to $d$, then
$$\rho(I^{\left[d+1\right]}) \leq \max \left\lbrace\rho(I), \rho(I^{\left[d\right]}), \dfrac{2}{3}(\rho(I)+\rho(I^{\left[d\right]})) \right\rbrace. $$
On the other hand
$$\rho(I^{\left[d\right]} < 2\rho(I),\;\; \dfrac{2}{3}(\rho(I)+\rho(I^{\left[d\right]})) < \dfrac{2}{3}(\rho(I)+2\rho(I))= 2\rho(I). $$
Hence $\rho(I^{\left[d+1\right]}) < 2\rho(I)$. The proof is concluded.
\end{proof}

\section{How large can $\rho(I)-\rho_a(I)$ be?}
\label{sect_difference}
For integers $m,d$ such that $1\le m\le d$, let $I_{m,d}$ be the ideal generated by products of $d-m+1$ different variables in the polynomial ring $A=\kk[x_1,\ldots,x_d]$. It is not hard to check that we have a primary decomposition
\begin{equation}
\label{eq_primdecomp}
I_{m,d}=\bigcap_{1\le i_1<\cdots <i_m\le d} (x_{i_1},\ldots,x_{i_m}).
\end{equation}
The following result due to Lampa-Baczy\'nska--Grzegorz determines the asymptotic resurgence of $I_{m,d}$. Villarreal has recently obtained a more complete statement in \cite[Theorem 3.19]{Vil2022}.
\begin{prop}[{\cite[Theorem C]{LBG15}}]
\label{prop_rho_a_uniform_md_ideal}
For all integers $m,d$ such that $1\le m\le d$, there is an equality
\[
\rho_a(I_{m,d})=\frac{m(d-m+1)}{d}.
\]
\end{prop}
\begin{cor}
\label{cor_equigen_ideal}
Let $m\ge 2$ be an integer, $A=\kk[x_1,\ldots,x_{2m-1}]$ and $I=I_{m,2m-1}$ be the ideal generated by all the products of $m$ distinct variables:
\[
I=(x_{i_1}x_{i_2}\cdots x_{i_m}: 1\le i_1 <i_2<\cdots<i_m \le 2m-1).
\]
Then $x_1\cdots x_{2m-1} \in I^{(m)}\setminus I^2$ and $\rho_a(I)=\dfrac{m^2}{2m-1}$.
\end{cor}
\begin{proof}
 It is clear that $x_1\cdots x_{2m-1} \notin I^2$. Thanks to \eqref{eq_primdecomp}, one has
\[
I^{(m)} = \bigcap_{1\le i_1<i_2<\cdots <i_m \le 2m-1} (x_{i_1}, x_{i_2},\ldots, x_{i_m})^m.
\]
Because $x_1\cdots x_{2m-1} \in (x_{i_1}, x_{i_2},\ldots, x_{i_m})^m$ for all $1\le i_1<i_2<\cdots <i_m \le 2m-1$, we obtain the first assertion.

That $\rho_a(I)= m^2/(2m-1)$ follows from Proposition \ref{prop_rho_a_uniform_md_ideal} with $d=2m-1$.
\end{proof}

The following result answers in the positive \cite[Conjecture 3.9]{BHJT}.
\begin{thm}
\label{thm_large_difference}
There exists a sequence of polynomial rings $R_n$ and squarefree monomial ideal $P_n\subseteq R_n$ generated by forms of the same degree, such that $\rho(P_n)-\rho_a(P_n)\to \infty$ when $n \to \infty$.
\end{thm}

\begin{proof}
For each integer $m\ge 2$, choose $I$ as in  \Cref{cor_equigen_ideal}. Using Lemma \ref{lem_sum_and_resurgence} for $v=2,m_1=m_2=m_3=m$, we get a squarefree monomial ideal $P_m$ with
\[
\rho(P_m)\ge \frac{3m}{4}, \rho_a(P_m)=\frac{m^2}{2m-1}.
\]
The equality follows from Equation \eqref{eq_rhoa_sum}. Hence $\frac{2m^2-3m}{4(2m-1)} \le \rho(P_m)-\rho_a(P_m) \to \infty$ when $m\to \infty$. Since $I$ is generated by forms of the same degree $m$, so is $P_m$.
\end{proof}

\section{On the possible values of $\rho(I+J)$ given $\rho(I)$ and $\rho(J)$}
\label{sect_resset}
Let $a,b\ge 1$ be real numbers, and
\begin{gather*}
\Res(a,b)=\left\{\rho(I+J):  \text{$I$ and $J$ are nonzero proper homogeneous ideals living in} \right.\\
                                    \text{disjoint sets of variables such that} \left. \rho(I)=a, \rho(J)=b \right\}. 
\end{gather*}
We have
\begin{prop}
Let $a,b\ge 1$ be real numbers. Denote $M=\max\{a,b\}, m=\min\{a,b\}$. The following statements hold.
\begin{enumerate}[\quad \rm (1)]
\item There is a containment $\Res(a,b)\subseteq \left[M, \max\left\{M, \frac{2(a+b)}{3}\right\}\right]$.
 \item If $M\ge 2m$ and $\Res(a,b)\neq \emptyset$ then $\Res(a,b)=\{M\}$.
 \item There is an equality 
 $$
 \Res(1,1)=\{1\} \bigcup \left\{\frac{n+1}{n}: n\ge 3 \right\}.
 $$
\end{enumerate}
\end{prop}
\begin{proof}
(1) and (2) are immediate from \Cref{thm_sharpbound_rho_sum}.

(3): By \cite[Theorems 3.6 and 3.9]{JKM2021}, the containment
\[
 \Res(1,1) \subseteq \{1\} \bigcup \left\{\frac{n+1}{n}: n\ge 3 \right\}
\]
holds. Note that the arguments for \emph{ibid.} use the ``associated'' binomial expansion for symbolic powers, which is by now validated by \Cref{thm.symbBIN}. That $1\in \Res(1,1)$ is immediate by choosing $I$ and $J$ to be isomorphic to $(x)\subseteq \kk[x]$. Let $I_d=(x^{2d+1},x^{2d-1}y^2,y^{2d+1}z)$ be the ideal of \Cref{lem_rho1_3generators}.  For an integer $n\ge 3$, choosing $I=I_1$ and $J \cong I_{n-2}$ so that they live in disjoint sets of variables. Applying \cite[Theorem 3.9]{JKM2021} for $r=2, p_1=2, p_2=n-1$, we get $\rho(I+J)=\dfrac{n+1}{n} \in \Res(1,1)$.
\end{proof}
It would be interesting to determine the set $\Res(a,b)$ in further cases, for example $\Res(n,n)$ where $n\ge 2$ is an integer. We wonder whether for any rational number $a\ge 1$, $\Res(a,a)$ is always a countable set.

\begin{rem}
The rationality of resurgence and asymptotic resurgence is an interesting topic considered in  \cite{DFMS, DD21, JKM2021}. We do not know whether if $a,b\in \Q$, then $\Res(a,b)\subseteq \Q$, namely if $\rho(I), \rho(J)\in \Q$ then so is $\rho(I+J)$. This has a positive answer if $\rho(I)=\rho(J)=1$ \cite[Theorem 3.9]{JKM2021} or if $I$ and $J$ are both monomial ideals \cite[Theorem 3.7]{DD21}, as the symbolic Rees algebras of monomial ideals are known to be noetherian. By \Cref{thm_sharpbound_rho_sum}, we also have a positive answer in the case $\max\{\rho(I),\rho(J)\} \ge 2\min\{\rho(I),\rho(J)\}$.
\end{rem}

\subsection*{Acknowledgments}
This work is partially supported by a grant from the International Centre for Research and Postgraduate Training in Mathematics, VAST (grant number ICRTM03$\_$2020.05). The second author acknowledges the generous support from the Simons Foundation Targeted Grant for the Institute of Mathematics - VAST (Award number: 558672), and from the Vietnam Academy of Science and Technology (grants CSCL01.01/22-23 and NCXS02.01/22-23).  He is also grateful to T\`ai H\`a and Arvind Kumar for inspiring discussions related to the content of this work. Finally, we are indebted to the information provided by Macaulay2  experiments \cite{GS96} at various stages of this project.


\begin{thebibliography}{99}

\bibitem{BHJT}
S. Bisui, H.T. H\`a, A.V. Jayanthan, and A.C. Thomas,
\emph{Resurgence numbers of fiber products of projective schemes}.
 Collect. Math. {\bf 72} (2021), no. 3, 605--614. 

	\bibitem{BH10}
C. Bocci and B. Harbourne,
\emph{Comparing powers and symbolic powers of ideals}.
 J. Algebraic Geom., {\bf 19} (2010), no. 3, 399--417.

 \bibitem{CHHVT}
E. Carlini, H.T. H\`a, B. Harbourne and A. Van Tuyl,
\emph{Ideals of powers and powers of ideals: Intersecting Algebra, Geometry and Combinatorics}.
Lecture Notes of the Unione Matematica Italiana, Vol. 27. Springer International Publishing, 2020.

\bibitem{DDS+}
H. Dao, A. De Stefani, E. Grifo, C. Huneke, and L. N\'u\~nez-Betancourt, 
\emph{Symbolic powers of ideals}.
In \emph{Singularities and Foliations. Geometry, Topology and Applications}, 387–432, Springer Proc. Math.
Stat., {\bf 222}, Springer, Cham, 2018.

\bibitem{DD21}
M. DiPasquale and B. Drabkin,
\emph{On resurgence via asymptotic resurgence}.
J. Algebra {\bf 587} (2021), 64--84.


\bibitem{DFMS}
M. DiPasquale, C.A. Francisco, J. Mermin, and J. Schweig,
\emph{Asymptotic resurgence via integral closures}.
Trans. Amer. Math. Soc. {\bf 372} (2019), no. 9, 6655--6676.

\bibitem{DHNSTT}
M. Dumnicki, B. Harbourne, U. Nagel, A. Seceleanu, T. Szembert and H. Tutaj-Gasi\'nska,
\emph{Resurgences for ideals of special point configurations in $\Pb^N$ coming from hyperplane arrangements}.
J. Alg. {\bf 443} (2015), 383--394.


\bibitem{ELS01}
L. Ein, R. Lazarsfeld, and K.E. Smith.
\emph{Uniform bounds and symbolic powers on smooth varietie}s.
Invent. Math. {\bf 144} (2001), 241--252.

\bibitem{GS96}
D. Grayson and M. Stillman, 
\emph{Macaulay2, a software system for research in algebraic geometry}. 
Available at \textsf{http://www.math.uiuc.edu/Macaulay2}. 

\bibitem{GHM2020}
E. Grifo, C. Huneke and V. Mukundan,
\emph{Expected resurgences and symbolic powers of ideals}.
J. Lond. Math. Soc. (2), {\bf 102} (2020), no. 2, 453--469.

\bibitem{GSV2022}
G. Grisalde, A. Seceleanu, and R.H. Villarreal,
\emph{Rees algebras of filtrations of covering polyhedra and integral closure of powers of monomial ideals}.
Res. Math. Sci. {\bf 9} (2022), Article number: 13. 

	\bibitem{GHV13}
E. Guardo, B. Harbourne, and A. Van~Tuyl.
\emph{Asymptotic resurgences for ideals of positive dimensional subschemes
of projective space}.
 Adv. Math. {\bf 246} (2013), 114--127.



\bibitem{HJKN}
H.T. H\`a, A.V. Jayanthan, A. Kumar, and H.D. Nguyen,
\emph{ Binomial expansion for saturated and symbolic powers of sums of ideals}, 
Preprint (2021), arXiv:2112.09338 [math.AC].

\bibitem{HNTT}
H.T. H\`a, H.D. Nguyen, N.V. Trung and T.N. Trung,
\emph{Symbolic powers of sums of ideals}.
Math. Z. {\bf 294} (2020), 1499--1520.

\bibitem{HTr19}
H.T. H\`a and N.V. Trung,
\emph{Membership criteria and containments of powers of monomial ideals}.
Acta Math. Vietnam. {\bf 44} (2019), 117--139 


\bibitem{HH02}
M. Hochster and C. Huneke.
\emph{Comparison of symbolic and ordinary powers of ideals}.
Invent. Math., {\bf 147} (2002), no. 2, 349--369.

\bibitem{JKM2021}
A.V. Jayanthan, A. Kumar and V. Mukundan,
\emph{On the resurgence and asymptotic resurgence of homogeneous ideals}.
Math. Z. (2022). \textsf{https://doi.org/10.1007/s00209-022-03138-w}.

\bibitem{LBG15}
M. Lampa-Baczy\'nska and M. Grzegorz,
\emph{On the containment hierarchy for simplicial ideals}.
J Pure Appl. Algbebra {\bf 219} (2015), no. 12, 5402--5412.

\bibitem{Vil2022}
R.H. Villarreal,
\emph{A duality theorem for the ic-resurgence of edge ideals}.
Preprint (2022), arXiv:2203.01268 [math.AC].

\end{thebibliography}
\end{document}